\documentclass{amsart}

\newtheorem{theorem}{Theorem}
\newtheorem{corollary}[theorem]{Corollary}
\newtheorem{thm}{Theorem}[section]
\newtheorem{cor}[thm]{Corollary}
\newtheorem{lem}[thm]{Lemma}
\newtheorem{prop}[thm]{Proposition}
\theoremstyle{definition}

\theoremstyle{remark}

\newtheorem{examp}{Example}
\numberwithin{equation}{section}

\newcommand{\eps}{\varepsilon}

\newcommand{\R}{\mathbb{R}}

\newcommand{\C}{\mathbb{C}}

\newcommand{\N}{{\mathbb N}}
\newcommand{\Z}{{\mathbb Z}}

\def \eqskip { \vspace*{2mm} }

\newcommand{\om}{{\omega}}

\newcommand{\De}{{\Delta}}

\newcommand{\la}{{\lambda}}

\newcommand{\atanh}{{\rm \ atanh }}

\newcommand{\disps}{\displaystyle}

\newcommand{\im}{{\sqrt{-1}}}

\newcommand{\ov}{\overline}
\newcommand{\op}{{\overline{\partial}}}
\newcommand{\p}{{\partial}}

\begin{document}

\title[$S^{1}-$invariant metrics on $S^{2}$]
{On the invariant spectrum of $S^{1}-$invariant
metrics on $S^{2}$}

\author{Miguel Abreu and Pedro Freitas}

\address{Departamento de Matem\'{a}tica, Instituto Superior T\'{e}cnico,
Av.Rovisco Pais, 1049-001 Lisboa, Portugal}
\email{mabreu@math.ist.utl.pt, pfreitas@math.ist.utl.pt}

\thanks{The authors were partially supported by
FCT (Portugal) through program POCTI, and grants
PCEX/C/MAT/44/96 and POCTI/32931/MAT/2000, respectively.}



\date{\today}

\dedicatory{}


\begin{abstract}
A theorem of J.\ Hersch (1970) states that for any smooth metric on $S^2$,
with total area equal to $4\pi$, the first nonzero eigenvalue of the
Laplace operator acting on functions  is
less than or equal to $2$ (this being the value for the standard round metric).
For metrics invariant under the standard $S^1$-action on $S^2$,
one can restrict the Laplace operator to the subspace of $S^1$-invariant functions
and consider its spectrum there. The corresponding
eigenvalues will be called invariant eigenvalues, and the purpose of 
this paper is to analyse their possible values.

We first show that there is no general analogue of Hersch's
theorem, by exhibiting explicit families of $S^1$-invariant metrics
with total area $4\pi$
where the first invariant eigenvalue ranges through any value between
$0$ and $\infty$. We then restrict ourselves to 
$S^1$-invariant metrics that can be embedded in ${\bf R}^3$ as surfaces of
revolution. For this subclass we are able to provide optimal upper bounds for all 
invariant eigenvalues. As a consequence, we obtain an analogue of Hersch's
theorem with an optimal upper bound (greater than $2$ and geometrically
interesting). This subclass of metrics on $S^2$ includes all $S^1$-invariant
metrics with non-negative Gauss curvature.

One of the key ideas in the proofs of these results comes from symplectic geometry,
and amounts to the use of the moment map of the $S^1$-action as a coordinate
function on $S^2$.
\end{abstract}

\maketitle

\section{Introduction}

Let $S^2$ be the $2$-sphere and $g$ any smooth Riemannian metric on
it. Denote by $\la(g)$ the first non-zero eigenvalue of the
Laplace operator defined by $g$, acting on functions on $S^2$. In
1970, J.Hersch~\cite{hers} showed that
$$\la(g) \leq \frac{8\pi}{\mbox{Area}_{g}(S^2)}\ .$$
If we scale the metric so that its total volume is the standard $4\pi$,
then Hersch's theorem becomes
\begin{equation}\label{eq:hersch}
\la(g)\leq 2\ .
\end{equation}
This was generalized to surfaces of higher genus by P.Yang and 
S.T.Yau~\cite{yaya} and to K\"ahler metrics on projective complex 
manifolds by J.-P.Bourguignon, P.Li and S.T.Yau~\cite{bly}

The purpose of this paper is to analyze when and how results of this 
type (i.e. upper bounds for eigenvalues)
can be obtained in an invariant setting. We consider only smooth metrics $g$
on $S^2$ with total area $4\pi$ and invariant under the standard
$S^1$-action (fixing the North and South poles). Denote by
$\la_{j}(g)$, with
$$0=\la_{0}(g) < \la_{1}(g) < \la_{2}(g) < \cdots\ ,$$
the invariant eigenvalues of the Laplace
operator defined by $g$ (i.e. with $S^1$-invariant eigenfunctions).
What can be said about the possible values of $\la_{j}(g)$?

Our first result says that there are no general restrictions on 
$\la_{1}(g)$ (and so, no general analogue of Hersch's theorem).

\begin{theorem}\label{thm1}
Within the class of smooth $S^1$-invariant metrics $g$ on $S^2$ with total
area $4\pi$, the first invariant eigenvalue $\la_{1}(g)$ can be any
number strictly between zero and $\infty$.

The same is true within the subclass of those metrics that have fixed
Gauss curvature at the poles: $K_{g}(N) = K_{g}(S) = 1$.
\end{theorem}

Our second result gives optimal upper bounds on the possible values 
of the invariant eigenvalues for the subclass of
$S^1$-invariant metrics on $S^2$ that are isometric to a surface
of revolution in $\R^3$. In particular, it includes an analogue of
Hersch's theorem within this geometrically interesting subclass of metrics.

\begin{theorem}\label{thm2}
Within the class of smooth $S^1$-invariant metrics $g$ on $S^2$ with total
area $4\pi$ and corresponding to a surface of revolution in
$\R^3$, we have that
$$\la_{j}(g) < \frac{\xi_{j}^2}{2},\ \ j=1,\ldots,$$
where $\xi_{j}$ is the $\left((j+1)/2\right)^{\rm th}$ positive zero of the
Bessel function $J_{0}$ if $j$ is odd, and the $\left(j/2\right)^{\rm th}$ 
positive zero of $J_{0}'$ if $j$ is even. These bounds are optimal.

In particular,
$$\la_{1}(g) < \frac{\xi_{1}^2}{2} \approx 2.89.$$
\end{theorem}

As a consequence of this last theorem, we have the following
\begin{corollary}
\label{cor1}
Let $g$  be any $S^1$-invariant metric on $S^2$ with total area $4\pi$
and non-negative Gauss curvature. Then the eigenvalues $\la_{j}$ 
satisfy the same upper bounds as in the previous theorem. These bounds 
are optimal.
\end{corollary}
\begin{proof}
Either apply Theorem~\ref{thm2}, using the fact that any
$S^1$-invariant metric on $S^2$ with non-negative scalar curvature is
isometric to a compact surface of revolution in $\R^3$
(see e.g.~\cite{kawa} and the references therein), or see at the end of \S 4 why
the proof of Theorem~\ref{thm2} also proves this corollary (without
using the above fact).
\end{proof}

Since the problem we are dealing with is essentially a 
$1$-dimensional eigenvalue problem, one might think that the size of 
the invariant eigenvalues would be related to the length of a geodesic 
meridian joining the North and South poles (which is the same as the 
{\it diameter} of $S^2$ with the given $S^1$-invariant metric):
{\it small diameter} $\leftrightarrow$ {\it large eigenvalues} and 
{\it large diameter} $\leftrightarrow$ {\it small eigenvalues}.
However we will see in \S 5 that this is not the case. Although the
problem is $1$-dimensional, it has curvature coming from the 
underlying $2$-sphere and that plays an important role.

The simplest family of $S^1$-invariant metrics for which 
$\la_{1}(g)\rightarrow\infty$, considered in the proof of
Theorem~\ref{thm1}, does have the property that the corresponding 
diameter $D(g)$ tends to zero (while the total volume is always 
$4\pi$ !). However, we will see that it is possible to 
change this particular family of metrics in two different ways:
\begin{itemize}
\item[(i)] in the first we keep the property that 
$D(g)\rightarrow 0$ but change the first invariant eigenvalue so 
that now $\la_{1}(g)$ also tends to zero;
\item[(ii)] in the second we keep the property that $\la_{1}(g)\rightarrow\infty$
but change the diameter so that now $D(g)$ also tends to $\infty$.
\end{itemize}

We will also see in the paper (\S 4) why the subclass of $S^1$-invariant 
metrics corresponding to surfaces of revolution in $\R^3$ is much 
more rigid than the full class of all $S^1$-invariant metrics on $S^2$.
We will quantify that rigidity very precisely, and conclude that among 
surfaces of revolution in $\R^3$ with total area $4\pi$, the one 
maximizing all invariant eigenvalues also minimizes the diameter and 
consists of the union of two flat discs of area $2\pi$ each (a 
singular surface). 

The above minimization property of the diameter of the double flat 
disc is just a generalization in this particular $S^1$-invariant 
setting of the following conjecture of Alexandrov: for every 
$(S^2,g)$ with non-negative curvature the ratio 
$\mbox{Area}_{g}(S^2)/D^{2}(g)$ is bounded from above by $\pi/2$,
and this value is only attained in the singular case of the double 
flat disc. 
According to~\cite{berg} this conjecture is still open. The $S^1$-invariant
case sugests that another possible interesting statement for it would 
be obtained by replacing {\it ``with non--negative curvature''} 
with {\it ``isometric to a closed surface in $\R^3$''} (recall that, 
by~\cite{nire}, any $(S^2,g)$ with positive curvature is isometric 
to a closed surface in $\R^3$).

It will also become clear in the paper (see in particular \S 4) that 
as one deforms the standard sphere towards the union of two 
flat discs, through a family of positive curvature ellipsoids of 
revolution intuitively obtained by ``pressing'' the North and South 
poles against each other, the first invariant eigenvalue increases
from $2$ to the limiting value $\frac{\xi_{1}^2}{2} \approx 2.89$.
Due to Hersch's Theorem and to the fact that the first noninvariant 
eigenvalue has multiplicity two, it follows that for small deformations of
this type the first invariant eigenvalue (which is larger than two) is the
third eigenvalue in the full spectrum. Since for the standard sphere the first
eigenvalue is 
equal to two and has multiplicity three, we obtain that any of these slightly 
deformed metrics have their third eigenvalue (in the full spectrum) 
larger than the third eigenvalue of the standard sphere. In this way,
we have obtained examples of metrics with positive curvature which 
provide a negative answer to the following question raised by S.T. Yau
(Problem \#71 in~\cite{yau}).

\ \newline 
{\em Let $0=\mu_{0}(g) < \mu_{1}(g) \leq \cdots \leq \mu_{m}(g) \leq 
\cdots $ be the (full) spectrum of any $(S^2,g)$. Is $\mu_{m}(g) \leq
\mu_{m}(\mbox{\rm standard}),\ \forall m$?}

\ \newline 
\noindent We note that another example providing a negative answer to this 
question had already been given by Engman in~\cite{eng1}, but in that case 
the metric had some negative curvature.

The proofs of Theorems~\ref{thm1} and~\ref{thm2} are an illustration of 
the usefulness of symplectic coordinates for some problems in Riemannian
geometry. In dimension $2$, any Riemannian metric determines a symplectic
form (just the corresponding area form), and symplectic coordinates are any
coordinates $(x,y)$ on which this symplectic form is the standard
$dx\wedge dy$. The existence of a circle action of isometries,
determines a particularly nice choice of symplectic coordinates (i.e.
action/angle coordinates) and these turn out to be very convenient for
the problem at hand.

This remark describes the simplest particular case of a more general
picture valid for any symplectic toric manifold (i.e. a symplectic
manifold of dimension $2n$ with an effective Hamiltonian action of the
$n$-torus $T^n$). The interested reader can look in~\cite{guil}, where
this general picture was used for the first time (to describe some K\"{a}hler
metrics on toric varieties), and~\cite{abr}.

In the $S^2$-setting, M.\ Engman has also used this type of 
coordinates to derive spectral properties of $S^1$-invariant metrics. 
In particular, he proves in~\cite{eng1} that the first invariant 
eigenvalue can be arbitrarily large (using a trace formula, while we 
use a Hardy type inequality to prove Theorem~\ref{thm1}), and 
in~\cite{eng2} he obtains the value $3$ as an upper bound for the 
first invariant eigenvalue of surfaces of revolution in $\R^3$. We 
thank Rafe Mazzeo for pointing out M.\ Engman's work to us.

The paper is organized as follows. In \S 2 we describe in detail the
particular coordinates on $S^2$ that will be used, and how all the
relevant Riemannian quantities can be expressed in this way.
Theorem~\ref{thm1} is proved in \S 3, while \S 4 is devoted to the
proof of Theorem~\ref{thm2}. Examples that show the general independence 
between the first invariant eigenvalue and the diameter are given in \S 5. 
In the Appendix we present the explicit solutions of a variational 
second order ODE that comes up naturally in \S 3.

\section{Preliminaries and background}

\subsection{Description of invariant metrics}

Let $S^2\subset\R^3$ be the standard sphere of radius $1$, with an
$S^1$-action given by rotation around the vertical axis. An
equivariant version of the Uniformization Theorem says that there is
only one $S^1$-invariant conformal structure on $S^2$. More
explicitly, given any $S^1$-invariant metric $g$ on $S^2$ there is an
$S^1$-equivariant diffeomorphism $\varphi:S^2\to S^2$ such that
$\varphi^\ast (g)$ is pointwise conformally equivalent to the standard
metric $g_0$ on $S^2$. This means that $\varphi^\ast (g)$ is
compatible with the standard ($S^1$-invariant) complex structure $j_{0}$
on $S^2$. Since we are in real dimension $2$, this is equivalent to
saying that $\varphi^\ast (g)$ is a K\"{a}hler metric on the complex
manifold $(S^2 , j_{0})$.

We conclude that in order to study spectral properties of
$S^1$-invariant metrics on $S^2$, it is enough to consider only
$S^1$-invariant K\"{a}hler metrics on the complex manifold $(S^2 ,
j_{0})$. Moreover, by multiplying a given metric by an appropriate
constant, we may consider only metrics with a fixed total volume (say
$4\pi$). We will now give an explicit description of this class of
metrics.

Let $g$ be an $S^1$-invariant K\"{a}hler metric on $(S^2 , j_{0})$,
with associated K\"{a}hler (or area) form denoted by $\omega$, and
such that $\mbox{Vol}_g (S^2) = \int_{S^2} \omega = 4\pi$. Because
$S^2$ is simply connected, the $S^1$-action is Hamiltonian with
respect to $\omega$ and we denote a corresponding Hamiltonian
function (or moment map) by $H:S^2\to\R$. A consequence of the
Duistermaat-Heckman theorem in this very simple setting is that the
push-forward by $H$ of the measure determined on $S^2$ by $\omega$
is the measure $\mu$ on $\R$ given by
$$ \mu(A) = 2\pi\cdot m(A\cap H(S^2)),\ A\subset\R,$$
where $m$ is standard Lebesgue measure. Hence we have that the image
interval $H(S^2)$ has length $2$. Since $H$ is only determined up
to a constant, we will normalize it so that $H(S^2)=[-1,1]\subset\R$
(this is equivalent to requiring that $\int_{S^2}H\cdot\omega = 0$).
This moment polytope $P = [-1,1]$ is determined by the affine
functions
$$l_1(x)=1+x\ \mbox{ and }\ l_2(x)=1-x$$
in the sense that $x\in P$ if and only if $l_i(x)\geq 0,\
i=1,2$, and $x\in P^{\circ}\equiv \mbox{interior of}\ P = (-1,1)$
if and only if $l_i(x)>0,\ i=1,2$. These two affine functions will be
relevant below.

The inverse image $H^{-1}(P^\circ)$ (i.e. the sphere minus the
two poles fixed by the $S^1$-action)
is the complex torus $M=\C/2\pi\im\Z$,
and $S^1 = \R/2\pi\Z$ acts on $M$ by the action:
$$
S^1\times M\rightarrow M,\ (t,z)\mapsto z + \im t
$$
where $z=u+\im v \in M$ and $t\in S^1$.
If $\om$ is an $S^1$-invariant K\"{a}hler form on $M$, for which the
$S^1$ action is Hamiltonian, then there exists a function $F=F(u),\ u={\rm Re}
\ z$, such that $\om = 2\im\p\op F$ and the moment map $H : M
\rightarrow\R$ is given by $H(z)=d F/d u$.

It follows that $\om$ can be written in the form
\begin{equation} \label{eq:Kahlerform}
\frac{\im}{2}\frac{d^2 F}{d u^2} dz\wedge d \overline{z}
\end{equation}
and the restriction to $\R  (={\rm Re}\ \C)$ of the K\"{a}hler
metric is the Riemannian metric
\begin{equation} \label{eq:Riemetric}
\frac{d^2 F}{d u^2}\ d u^2\ .
\end{equation}
Under the Legendre transform given by the moment map
\begin{equation} \label{eq:legendrex}
x=\frac{d F}{d u} = H\ ,
\end{equation}
this is the pull-back of the metric
$$
\frac{d^2 G}{d x^2}\ d x^2
$$
on $P^\circ$, where the potential $G$ is the Legendre function
dual to $F$ (up to a linear factor in $x$). More explicitly, the
inverse of the Legendre transform (\ref{eq:legendrex}) is
\begin{equation} \label{eq:legendreu}
u = \frac{d G}{d x} + a,
\end{equation}
with $a\in\R$ a constant. It follows from (\ref{eq:legendrex}) and
(\ref{eq:legendreu}) that
\begin{equation} \label{eq:relxu}
\frac{d^2 G}{d x^2} = \left(\frac{d^2 F}{d u^2}\right)^{-1}\
\mbox{at}\  x = \frac{d F}{d u}\ .
\end{equation}

For the standard metric $g_0$ on $S^2$, the function $G_0$ obtained
in this way is (see~\cite{guil})
$$ G_0 = \frac{1}{2}\sum_{k=1}^{2} l_k(x)\log l_k(x)
       = \frac{1}{2}[(1+x)\log (1+x) + (1-x)\log (1-x)] $$
and the metric induced on $P^{\circ}$ is
$$g_0(x)\ \mbox{d}x^2 \ \ \mbox{ with }\ \ g_0(x) = \frac{d^2 G_0}{d x^2}
= \frac{1}{1-x^2}\ .$$

For any other $S^1$-invariant K\"{a}hler metric $g$ on $S^2$, with
total volume $4\pi$, we have that the corresponding volume form
$\omega$ satisfies
$$\omega - \omega_0 = \frac{\im}{2}\p\op\phi\ , $$
where $\omega_0$ denotes the volume form of the standard metric
$g_0$ and $\phi$ is a smooth real-valued $S^1$-invariant function
on $S^2$. It follows from the above construction that the
associated potential $G$ is given by $$ G = G_0 + \Phi\ ,$$ with
$\Phi \in C^{\infty}(P)$, and hence the metric induced by $g$
on $P^{\circ}$ is of the form
$$g(x)\ \mbox{d}x^2 \ \ \mbox{ with }\ \ g(x) = \frac{d^2 G}{d x^2}
= \frac{1}{1-x^2} + h(x)\ ,$$
where $h\in C^{\infty}(P)$ is such that $g(x)>0$ for $x\in P^{\circ}$.

We conclude from this construction that the space of all
$S^1$-invariant K\"{a}hler metrics $g$ on $S^2$ with total volume
$4\pi$  may be identified with the space of all functions
\begin{equation} \label{eq:defh}
h\in C^{\infty}(P)\ \mbox{such that} \ h(x)>-\frac{1}{1-x^2},\
x\in P^{\circ}\ .
\end{equation}

\subsection{The Laplacian}

If $f=f(u)$ is an $S^1$-invariant function on
$M=H^{-1}(P^{\circ})$, its Laplacian with respect to the
K\"{a}hler metric defined by (\ref{eq:Kahlerform}) is given by
\begin{equation} \label{eq:Lapu}
\De f = - \left(\frac{d^2 F}{d u^2}\right)^{-1} \frac{d^2 f}{d u^2}
\end{equation}
(note that this is not the same as the Laplacian of the restriction
of $f$ to $\R (={\rm Re}\ \C)$ with respect to the Riemannian
metric defined by (\ref{eq:Riemetric})). To write this operator in
terms of the moment map coordinate $x$, we note that it follows
from (\ref{eq:legendrex}) and (\ref{eq:relxu}) that
\begin{equation} \label{eq:dudx}
\frac{d}{d u} = \frac{d x}{d u}\frac{d}{d x} =
\frac{d^2 F}{d u^2}\frac{d}{d x} = (\frac{d^2 G}{d x^2})^{-1}
\frac{d}{d x} = \frac{1}{g(x)}\frac{d}{d x}
\end{equation}
at $x = \frac{d F}{d u}$. Hence the Laplacian is given in the $x$
coordinate by
\begin{equation} \label{eq:Lapx}
\De f = - \frac{d}{d x}\left(\frac{1}{g}\frac{d f}{d x}\right) =
- \left(\frac{1}{g} f'\right)'\ .
\end{equation}

\subsection{The invariant eigenvalues}

Given any $S^1$-invariant K\"ahler metric $g$ on $(S^2,j_0)$, with
volume form $\om$ and total volume $4\pi$, and a function
$f\in C^\infty(S^2)$, let
\begin{equation} \label{eq:norm2}
\|f\|_{g}^2 = \int_{S^2} f^2(x)\,\om + \int_{S^2}
|\nabla f|^2_{g}\,\om\ .
\end{equation}
The completion of $C^\infty(S^2)$ with respect to the above norm is
the Sobolev space $L^2_{1}(S^2)$. The Laplacian $\De_{g}$ is a
self-adjoint elliptic operator on $L^2_{1}(S^2)$ with a discrete
non-negative spectrum.

Because the push-forward (by the moment map $H$) of the measure
determined on $S^2$ by $\om$ is simply $2\pi$ times the Lebesgue measure
on the polytope $P$, we have that for $S^1$-invariant functions
$f\in C^\infty(S^2)$ the above norm can be written in terms of the $x$
coordinate as
\begin{equation} \label{eq:norm1}
\|f\|_{g}^2 = 2\pi\left(\int_{-1}^{1} f^2(x)\,dx + \int_{-1}^{1}
\frac{(f'(x))^2}{g(x)}\,dx \right)\ ,
\end{equation}
where $g(x)\,dx^2$ is the metric induced by $g$ on $P^{\circ}$.
As in any Riemannian manifold, different metrics induce equivalent norms
and we denote by $X$ the completion of $C^{\infty}(P)$ with respect to
any of them. Then the Laplacian defined by~(\ref{eq:Lapx}) is a
self-adjoint operator on $(X, \|\cdot\|_{g})$,
having discrete non-negative spectrum consisting exactly of the
eigenvalues $\la_{j}(g)$ of $\De_{g}$ on $L^2_{1}(S^2)$ with $S^1$-invariant
eigenfunctions $f_{g,j}$ (i.e. the $S^1$-invariant spectrum).

By the Min-Max principle we then have that the nontrivial invariant
eigenvalues $0<\la_{1}<\la_{2}<\ldots,$ are given by
\begin{equation} \label{eq:quotient}
\lambda_{j} = \lambda_{j}(g) = {\disps \inf_{f\in X_{g,j}} }
\frac{\disps\int_{-1}^{1} \frac{\disps (f'(x))^2}{\disps g(x)}\,dx}
{\disps \int_{-1}^{1} f^{2}(x)\,dx}, \ \  \ j=1,\ldots,
\end{equation}
where $X_{g,j} = \{f\in X: \int_{-1}^{1}f(x)f_{g,k}\,dx=0, \ 
k=0,\ldots, j-1, \mbox{ and } f\neq
0\}$ ($f_{g,0}\equiv 1$). A nice feature of the use of the moment map
coordinate $x$ is that $X_{1}\equiv X_{g,1}$ and the denominator of the 
above quotient do not depend
on the metric $g$, since the relevant measure on the polytope $P$ is
always $2\pi$ times Lebesgue measure.

\subsection{The diameter}

The diameter $D(g)$ of $(S^2,g)$ is equal to the length of any geodesic 
meridian joining the North and South poles. This can be computed as 
the length of the polytope $P$ with respect to the metric induced by 
$g$, and hence is given by
\begin{equation}\label{eq:diameter}
D(g) = \int_{-1}^{1}\sqrt{g(x)}\,dx\ .
\end{equation}

\subsection{The Gauss curvature}

The Gauss curvature $K$ of the K\"ahler metric defined
by~(\ref{eq:Kahlerform}) is given by
$$K = - \frac{1}{2}\left(\frac{d^2 F}{d u^2}\right)^{-1}
\frac{d^2 \log(d^2 f/du^2)}{d u^2}\ ,$$
which, using~(\ref{eq:relxu}) and~(\ref{eq:dudx}), can be written in
terms of the moment map coordinate $x=dF/du$ as
\[
\begin{array}{lcl}
K & = & -\frac{\disps 1}{\disps 2}g\frac{\disps 1}{\disps g}
\frac{\disps d}{\disps dx}\left[
\frac{\disps 1}{\disps g}\frac{\disps d}{\disps dx}\log(1/g)\right]
\eqskip\\
 & = & -\frac{\disps 1}{\disps 2}
 \left(\frac{\disps 1}{\disps g}\right)''\ .
\end{array}
\]
Hence, non-negative Gauss curvature amounts to $(1/g)''\leq 0$.

\subsection{Example}

The standard metric $g_{0}$ on $S^2$ is given on $P^{\circ}$ by
$$ g_{0}=\frac{dx^2}{1-x^2}\ .$$
The corresponding Laplacian on $S^1$-invariant functions is
$$\De f = - ((1-x^2)f')'\ \ \mbox{for}\ \ f\in C^{\infty}(P)\ ,$$
with invariant eigenfunctions the well-known Legendre polynomials
and invariant spectrum $\la_{n} = n(n+1),\,n\in\N$. In particular
this implies that
\begin{equation} \label{ineq:std}
\frac{\disps \int_{-1}^{1} (1-x^{2})
\left[f'(x)\right]^{2}dx}{\disps \int_{-1}^{1} f^{2}(x)dx}
\geq \la_{1}(g_{0}) = 2 \ \ \mbox{for any}\ \ f\in X_{1}\ ,
\end{equation}
an inequality that will be used later on.

The diameter and Gauss curvature are given by
$$ D_{0} = \int_{-1}^{1}\frac{1}{\sqrt{1-x^2}}\,dx = \pi \ \ \mbox{and}\
\ K_{0}(x) = -\frac{1}{2}(1-x^2)'' =  1$$
as expected.

\section{Proof of Theorem~\ref{thm1}}

We saw in the previous section that any smooth $S^1$-invariant metric
on $S^2$ is determined by a positive function $g\in C^{\infty}
(P^\circ)$, $P^\circ = (-1,1)$, of the form
\begin{equation} \label{eq:g}
g(x) = \frac{\disps 1}{\disps 1 - x^2} + h(x)
\end{equation}
with $h\in C^{\infty}(P)$, $P = [-1,1]$. We also saw that the relevant
Riemannian information (Laplacian, $\la_{1}$, scalar curvature)
can be written explicitly in terms
of $\ov{g} = 1/g$. From~(\ref{eq:g}) we get that $\ov{g}$ is of the
form
\begin{equation} \label{eq:inverseg}
\ov{g}(x) = \frac{\disps 1}{\disps g(x)} = (1-x^2)[1 +
(1-x^2)\ov{h}(x)]
\end{equation}
where $\ov{h}\in C^{\infty}(P)$ is such that $\ov{g}(x)>0,\
x\in P^\circ$. Note that $\ov{g}\in C^{\infty}(P)$. The functions $h$ and
$\ov{h}$ are related  to each other by
$$ \ov{h}(x) = - \frac{\disps h(x)}{\disps 1 + (1-x^2) h(x)}\
\mbox{and}\ h(x) = - \frac{\disps \ov{h}(x)}{\disps 1 + (1-x^2)
\ov{h}(x)} $$
and they both satisfy the inequality
$$ h(x),\ov{h}(x) > - \frac{\disps 1}{\disps 1-x^2},\
x\in P^{\circ}\ .$$
Note that any function $\ov{g}\in C^{\infty}(P)$ of the form given
by~(\ref{eq:inverseg}) satisfies
\begin{equation} \label{eq:bcondg}
\ov{g}(-1)=0=\ov{g}(1) \ \mbox{and}\ \ov{g}'(-1)=2=-\ov{g}'(1)\ .
\end{equation}

\subsection{Large first invariant eigenvalue\label{larginv}}

It follows from~(\ref{eq:quotient}) that in order to make $\lambda_{1}$
large one should choose the function $g$ small, and hence $\ov{g}$ as 
large as possible. It is then clear from~(\ref{eq:inverseg}) that the
simplest way to achieve that is to choose
$$ \ov{h}(x) = \mu = \ \mbox{constant}\ > 0 $$
and analyze what happens to $\lambda_{1}$ as $\mu$ tends to $\infty$.

Hence we will consider the family of metrics
\[
\begin{array}{lllll}
g_{\mu}(x) & = & \frac{\disps 1}{\disps 1-x^{2}}
- \frac{\disps \mu}{\disps 1 + (1-x^{2})\mu}
 & = & \frac{\disps 1}{\disps (1-x^{2})(1+(1-x^{2})\mu)}>0,
\end{array}
\]
with $0<\mu\in\R$, for which we have that
\[
\begin{array}{lcl}
\lambda_{1}(g_{\mu}) & = & {\disps \inf_{f\in X_{1}} }\frac{\disps
\int_{-1}^{1} \frac{\disps 1}{\disps g_{\mu}(x)}
\left[f'(x)\right]^{2}dx}{\disps \int_{-1}^{1} f^{2}(x)dx}\eqskip\\
& = & {\disps \inf_{f\in X_{1}} }\left\{\mu \frac{\disps
\int_{-1}^{1} (1-x^{2})^{2}
\left[f'(x)\right]^{2}dx}{\disps \int_{-1}^{1} f^{2}(x)dx}+
\frac{\disps \int_{-1}^{1} (1-x^{2})
\left[f'(x)\right]^{2}dx}{\disps \int_{-1}^{1} f^{2}(x)dx}\right\}
\eqskip\\
 & \geq & \mu {\disps
\inf_{f\in X_{1}} }
\frac{\disps \int_{-1}^{1} (1-x^{2})^{2}
\left[f'(x)\right]^{2}dx}{\disps \int_{-1}^{1} f^{2}(x)dx} + 2,
\end{array}
\]
this last inequality being valid because of~(\ref{ineq:std}).
Thus, to prove that $\lambda_{1}(g_{\mu})\rightarrow\infty$
as $\mu\rightarrow\infty$, it is enough to show that this last
quotient is bounded away from zero. In order to do this, we shall
consider the inequality
\begin{equation}
\label{mainineq0}
\int_{0}^{1} (1-x^{2})^{2p} [f'(x)]^{2}dx\geq C\int_{0}^{1}f(x)^{2}dx.
\end{equation}
>From results in~\cite{opku} it follows that it holds if and only if
$p$ is less than or equal to $1$. Their proof is based on Minkowski's
and H\"{o}lder's inequalities. Here we follow a different approach
which allows us to obtain the optimal constant in the case where $p$
is one. We also recover the optimal constant when $p$ is $1/2$, which
corresponds to Legendre's equation.

\begin{lem} \label{lemineq}
Let $f\in X_{1}$ satisfy $f(0)=0$. Then, for all $p\in[0,1]$ there
exists $C=C(p)\geq1$ such that inequality~(\ref{mainineq0}) is
satisfied.

Furthermore, $C(1/2)=2$ (Legendre's equation) and $C(1)=1$ are optimal
constants in~(\ref{mainineq0}). For all other $p\in[0,1)$ we have that
the optimal constant $C_{opt}$ satisfies $C_{opt}(p)\geq m$, where $m$
is the infimum of the function $F:(0,1)\to\R$ defined by
\[
F(x) = 2 + \frac{\disps 1-x^{2}}{\disps x^{2}}\left[1-
(1-x^{2})^{1-2p}\right].
\]
\end{lem}
\begin{proof} From
\[
0 \leq \int_{0}^{1} \left[ -\frac{\disps (1-x^{2})^{1-p}}{\disps x}
f(x) + (1-x^{2})^{p}f'(x)\right]^2dx
\]
and integrating by parts we obtain that
\[
\int_{0}^{1} (1-x^{2})^{2p}[f'(x)]^{2}dx \geq \int_{0}^{1} F(x)
f^{2}(x)dx\geq m \int_{0}^{1} f^{2}(x)dx.
\]

In the case where $p=1/2$, we get $F(x)\equiv 2$. Since for $f(x)=x$
we have equality, it follows that $C_{opt}(1/2) = 2$. When $p=1$, we
have that $F(x)\equiv 1$. The solution of the corresponding
Euler--Lagrange equations which is given in the Appendix suggests that
we now consider the sequence of functions defined by
\[
f_{\eps}(x) = \frac{\disps x}{\disps \sqrt{1+\eps-x^2} }
\]
to get
\[
\int_{0}^{1} (1-x^{2})^2 [f_{\eps}'(x)]^{2}dx = \frac{\disps
1+\eps}{4}- \frac{\disps 8+\eps}{\disps 8(1+\eps)^2}+ \frac{\disps
8+8\eps+3\eps^2} {\disps 8(1+\eps)^{5/2}}\atanh\left(\frac{\disps
1}{\disps
\sqrt{1+\eps}} \right)
\]
and
\[
\int_{0}^{1} f_{\eps}^{2}(x)dx = -1 + \sqrt{1+\eps}\atanh\left(
\frac{\disps 1}{\disps \sqrt{1+\eps}}\right).
\]
>From this it follows that
\[
\lim_{\eps\to 0^{+}} \frac{\disps \int_{0}^{1} (1-x^{2})^2
[f_{\eps}'(x)]^{2}dx}{\disps \int_{0}^{1} f_{\eps}^{2}(x)dx} = 1
\]
and so $C_{opt}(1) = 1$.
\end{proof}

\begin{thm} Let $f\in X_{1}$.
Then
\begin{equation}
\label{mainineq1}
\int_{-1}^{1} (1-x^{2})^{2p} [f'(x)]^{2}dx\geq C
\int_{-1}^{1}f(x)^{2}dx,
\end{equation}
where $C=C(p)$ is as in Lemma~\ref{lemineq}. In particular,
$C(1)=1$ is optimal.
\end{thm}
\begin{proof} Assume that $f$ has zero average. Then
\[
\begin{array}{lll}
\int_{-1}^{1} (1-x^{2})^{2p} [f'(x)]^{2}dx & = &\int_{-1}^{1}
(1-x^{2})^{2p} \left\{ [f(x)-f(0)]'\right\}^{2}dx\eqskip \\ & \geq &
C\int_{-1}^{1}\left[f(x)-f(0)\right]^{2}dx,
\end{array}
\]
where the inequality follows from Lemma~\ref{lemineq}. We thus have
that
\[
\begin{array}{lll}
\int_{-1}^{1} (1-x^{2})^{2p} [f'(x)]^{2}dx & \geq & C\int_{-1}^{1}
f^2(x)dx - 2C f(0)\int_{-1}^{1}f(x)dx + 2C f^{2}(0)\eqskip \\ &
\geq & C\int_{-1}^{1} f^{2}(x)dx.
\end{array}
\]
\end{proof}

These results enable us to conclude that for the family of metrics
considered above, the first invariant eigenvalue satisfies
$\lambda_{1}=\lambda_{1}(g_{\mu})\geq \mu+2\,\rightarrow\,\infty$
as $\mu\rightarrow\infty$.

The Gauss curvature $K_{\mu}$ of this family of metrics is given
by (see \S 2.5)
$$ K_{\mu}(x) = -\frac{1}{2}[(1-x^2)(1+\mu (1-x^2))]'' =
1 + 2\mu (1-3x^2)\ . $$

In particular, at the poles fixed by the $S^1$-action
we have that the curvature blows up:

$$ K_{\mu}(1) = K_{\mu}(-1) = 1 - 4\mu \rightarrow -\infty\ \
\mbox{as}\ \ \mu\rightarrow \infty\ .$$
It is however possible to make the first invariant eigenvalue
arbitrarily large while keeping the curvature at the poles fixed.
Consider the family of metrics defined by
$$\ov{g}_{\rho}(x) = (1-x^2)[1 + \rho (1-x^2)^2]\ ,$$
and let us analyze what happens when $\rho \rightarrow\infty$. The
Gauss curvature at the poles is fixed (equal to $1$) while the first
invariant eigenvalue is given by
$$\lambda_{1}(g_{\rho}) =
{\disps \inf_{f\in X_{1}} }\left\{
\frac{\disps \int_{-1}^{1} (1-x^{2})\left[f'(x)\right]^{2}dx}
{\disps \int_{-1}^{1} f^{2}(x)dx} +
\rho \frac{\disps \int_{-1}^{1} (1-x^{2})^{3} \left[f'(x)\right]^{2}dx}
{\disps \int_{-1}^{1} f^{2}(x)dx}\right\}\ .$$ Although we now have
that
$${\disps \inf_{f\in X_{1}} } \frac{\disps \int_{-1}^{1} (1-x^{2})^{3}
\left[f'(x)\right]^{2}dx} {\disps \int_{-1}^{1} f^{2}(x)dx} = 0\ ,$$
we can use Cauchy-Schwartz to obtain
$$\left(\int_{-1}^{1} (1-x^{2})\left[f'(x)\right]^{2}dx\right) \cdot
  \left(\int_{-1}^{1} (1-x^{2})^{3} \left[f'(x)\right]^{2}dx\right)
  \geq\left(\int_{-1}^{1} (1-x^{2})^{2}
  \left[f'(x)\right]^{2}dx\right)^2 $$ which, together with~(\ref{ineq:std})
and Lemma~\ref{lemineq}, easily implies that
$$\la_{1}(g_{\rho})\geq\sqrt{4 + 2\rho} \rightarrow \infty\ \
\mbox{as}\ \ \rho\rightarrow\infty\ .$$

\subsection{Small first invariant eigenvalue}

We shall now see that it is possible to choose the metric in such a
way that the first invariant eigenvalue becomes arbitrarily small.
This can be achieved with a family of metrics ``dual'' to the previous
one. We choose
$$ h(x) \equiv \nu >0$$
and consider the family of metrics defined by
\[
\begin{array}{lllll}
g_{\nu}(x) & = & \frac{\disps 1}{\disps 1-x^{2}}+\nu,
\end{array}
\]
for constant $\nu$. The first invariant
eigenvalue of the Laplacian on $S^{2}$ corresponding to this family is
given by
\[
\begin{array}{lcl}
\lambda_{1}(g_{\nu}) & = & {\disps \inf_{f\in X_{1}} }\frac{\disps
\int_{-1}^{1} \frac{\disps 1-x^{2}}{\disps 1+\nu(1-x^{2})}
\left[f'(x)\right]^{2}dx}{\disps \int_{-1}^{1} f^{2}(x)dx}\eqskip\\
& < & \frac{\disps 1}{\disps \nu}{\disps \inf_{f\in X_{1}} }
\frac{\disps\int_{-1}^{1}
\left[f'(x)\right]^{2}dx}{\disps \int_{-1}^{1} f^{2}(x)dx}=
\frac{\disps \pi^{2}}{\disps 4\nu},
\end{array}
\]
and so we get that $\lambda_{1}(g_{\nu})\rightarrow 0$ as
$\nu\rightarrow\infty$.

A calculation shows that the Gauss curvature $K_{\nu}$ of this family 
of metrics is positive everywhere, tends to zero (as 
$\nu\rightarrow\infty$) for every $x\in(-1,1)$, while at the poles we 
have that
$$ K_{\nu}(1) = K_{\nu}(-1) = 1 + 4\nu \rightarrow \infty\ \
\mbox{as}\ \ \nu\rightarrow \infty\ .$$
The fact that one can get $\la_{1}$ arbitrarily small
with a family of metrics with fixed curvature at the poles
is left as an easy exercise to the reader. 

\section{Proof of Theorem~\ref{thm2}}

We now want to consider $S^1$-invariant metrics $g$ on $S^2$,
that correspond to closed surfaces of revolution in $\R^3$.
Such a surface is obtained by revolving a profile curve
$$ t\mapsto (0,p(t),q(t)),\ \mbox{for}\ 0\leq t \leq\ell,$$
about the third coordinate axis. We have necessarily that
$p(0)=0=p(\ell)$ and $p(t)>0$ for $0<t<\ell$. If the curve
is parametrized by arclength, which we will assume, then we
have in addition $\dot{p}(0) = 1 = - \dot{p}(\ell)$ and
\begin{equation} \label{eq:boundp}
(\dot{p}(t))^2 + (\dot{q}(t))^2 = 1,\ \mbox{for all}\ 0\leq t\leq\ell.
\end{equation}

Assuming that the total volume of such a surface is $4\pi$, then
the normalized moment map $H$ of the $S^1$-action, with respect
to the induced area form, ``projects'' the surface to the moment
polytope $P = [-1,1]\subset\R$. This can be seen as the
relation between the arclength coordinate $t$ and the moment map
coordinate $x$:
\begin{equation} \label{eq:relxt}
x = x(t) = H(t),\ t\in [0,\ell],\ x\in P\ .
\end{equation}
Under the above relation, the metric induced by $\R^3$ on the
surface of revolution gives rise to a metric $g$ on
$P^{\circ} = (-1,1)$,
$$g = g(x)\, dx^2$$
of the form considered in \S 2.1 . This function $g(x)$ determines
the inverse of~(\ref{eq:relxt}) by
\begin{equation} \label{eq:reltx}
t = t(x) = \int_{-1}^{x}\sqrt{g(s)}\,ds\ .
\end{equation}

\begin{prop} \label{prop:relpg}
The relation between $\ov{g} = 1/g$ and $p$ is given by
$$\ov{g}(x) = p (t(x))^2,\ x\in P\ . $$
\end{prop}
\begin{proof}
The curvature $K$ of the surface of revolution is given in the
arclength coordinate $t$ by~\cite{onei}
$$K(t) = - \frac{\ddot{p}(t)}{p(t)}\ ,$$
while in the moment map coordinate $x$ we have that (see \S 2.5)
$$ K(x) = - \frac{1}{2}\ov{g}''(x)\ . $$
This means that
\begin{equation} \label{eq:relpg1}
2 \ddot{p}(t(x)) = p(t(x))\cdot \ov{g}''(x)\ .
\end{equation}
Defining a function $\ov{p}:P\rightarrow \R$ by
$\ov{p}(x)=p(t(x))$, we get from~(\ref{eq:reltx})
and~(\ref{eq:relpg1}) the following differential equation relating
$\ov{p}$ and $\ov{g}$:
\begin{equation} \label{eq:relpg2}
\ov{g}\ov{p}'' + \ov{g}' \ov{p}' =
\frac{1}{2} (\ov{p}\ov{g}'' + \ov{g}' \ov{p}')\ .
\end{equation}
Given $\ov{p}$, one easily checks that the unique solution $\ov{g}$
of~(\ref{eq:relpg2}) satisfying
$$ \ov{g}(-1) = 0 = \ov{g}(1)\ \mbox{and}\ \ov{g}'(-1) = 2 =
-\ov{g}'(1) $$
is
$$ \ov{g} = \ov{p}^2$$
as required.
\end{proof}
\begin{cor} \label{cor:boundg}
For a surface of revolution in $\R^3$ with total area $4\pi$,
the corresponding metric on the moment polytope $P$,
$$g = g(x)\, dx^2 = \frac{1}{\ov{g}(x)}\, dx^2\ ,$$
is such that
$$ |\ov{g}'(x)| \leq 2,\ \forall x\in P\ .$$
\end{cor}
\begin{proof} From $\ov{g} = \ov{p}^2$ we have that
$$ \ov{g}' = 2\ov{p}\ov{p}'= 2\ov{p}(p\circ t)' = 2\ov{p}
(\dot{p}\circ t) t' = 2\ov{p}(\dot{p}\circ t)\frac{1}{\sqrt{\ov{g}}}
= 2(\dot{p}\circ t)\ . $$
Hence, using~(\ref{eq:boundp}), we get
$$ |\ov{g}'(x)| = 2 |\dot{p}(t(x))| \leq 2,\ \forall x\in P\ .$$
\end{proof}

We now know that, for a smooth closed surface of revolution in $\R^3$ with
total area $4\pi$, the corresponding function $\ov{g}\in
C^{\infty}(P)$ satisfies
\begin{equation} \label{embcond}
\ov{g}(-1)=0=\ov{g}(1),\ \ov{g}'(-1) = 2 = -\ov{g}'(1)\ \ \mbox{and}\ \
|\ov{g}'(x)|\leq 2,\ \forall x\in P\ .
\end{equation}
Any such $\ov{g}$ is clearly less than or equal to the ``tent'' function
$$\ov{g}_{\rm max}(x) = 2(1-|x|).$$ 
We shall now consider the invariant spectrum 
of this limit problem. The values of $\la_{j}(g_{\rm max})$ can
be explicitly determined in two different ways: one more geometric, the
other more analytic.

In the more geometric way, one interprets $g_{\rm max}$ as being the
singular metric on $S^2$ consisting of two flat discs, each of area
$2\pi$, glued along a singular equator. In fact, the scalar curvature
of $g_{\rm max}$ is zero where defined (i.e. $P\setminus\{ 0 \}$)
because $\ov{g}_{\rm max}$ is a polynomial of degree $1$ in $[-1,0)$
and $(0,1]$. Moreover, one can see $g_{\rm max}$ as the limit of the
family of surfaces of revolution in $\R^3$ consisting of ellipsoids of
total area $4\pi$ squeezed between two horizontal planes moving
towards each other (the associated family of functions $\ov{g}$ 
interpolates increasingly from $\ov{g}_{0}$ to $\ov{g}_{\rm max}$).
Because this is a family of ``even'' metrics (i.e. invariant under
the antipodal map or, equivalently, the corresponding functions $g$
and $\ov{g}$ are even on $P$), the associated invariant
eigenfunctions alternate between being odd and even, 
and hence either themselves or their first derivative
vanish at $x=0$ (the equator of
the surface of revolution). In the limit, and restricting attention
to just one of the hemispheres, the Laplace operator of this family
of metrics tends to the Euclidean Laplacian on a disc of area $2\pi$,
while the invariant eigenfunctions tend to its radially symmetric
eigenfunctions with either Dirichlet or Neumann boundary conditions. 
Hence we have that
$$\la_{j}(g_{\rm max}) = \frac{\xi_{j}^2}{2},\ \ j=1,\ldots,$$
where $\xi_{j}$ is the $\left((j+1)/2\right)^{\rm th}$ positive zero of the
Bessel function $J_{0}$ if $j$ is odd, and the $\left(j/2\right)^{\rm th}$ 
positive zero of $J_{0}'$ if $j$ is even. 
In particular,
$$\la_{1}(g_{\rm max}) = \frac{\xi_{1}^2}{2} \approx 2.89.$$

To obtain the same result from a more analytic perspective, we
consider the differential equation directly, that is,
\begin{equation} \label{eq:max}
\left\{
\begin{array}{ll}
[2(1-x)f']' + \lambda f = 0,\  x\in[0,1]\\ f(0) = 0 \mbox{ or }
f'(0) = 0,
\end{array}
\right.
\end{equation}
where the extra boundary condition in $f$ stems from the fact mentioned
above that we are considering an even metric and thus the eigenfunctions
alternate between being odd and even.

We now consider the change of variables defined by $t=\sqrt{2(1-x)}$,
which transforms this to the equation corresponding to the standard
Dirichlet/Neumann problem on the disk of radius $\sqrt{2}$. We are thus led to
the equation
\[
t^{2}\ddot{f}+t\dot{f} + \lambda t^{2}f = 0,
\]
with the conditions $f(\sqrt{2})=0$ or $f'(\sqrt{2})=0$. This equation
has $f(t)=J_{0}(\sqrt{\lambda} t)$ as a solution and we thus obtain the
above mentioned values for $\la_{j}(g_{\rm max})$.

The proof of Theorem~\ref{thm2} is now a direct consequence of 
the monotonicity principle (derived from Poincar\'{e}'s 
principle). Given any closed surface of revolution in $\R^3$
with total area $4\pi$, the corresponding function $\ov{g}\in C^{\infty}(P)$
satisfies~(\ref{embcond}), and hence $\ov{g}\leq \ov{g}_{\rm max}$
on $P$. This means that
$$\frac{\disps \int_{-1}^{1} 
\ov{g}(x) (f'(x))^2 dx} {\disps \int_{-1}^{1} f^{2}(x)dx}
\leq  \frac{\disps \int_{-1}^{1} 2 (1-|x|)
\left[f'(x)\right]^{2}dx} {\disps \int_{-1}^{1} f^{2}(x)dx}\ 
\mbox{for any}\ f\in X\ ,$$
which by the monotonicity principle implies the eigenvalue 
inequalities
\[
\la_{j}(g) \leq \la_{j}(g_{\rm max}).
\]
Moreover, because for a smooth surface of revolution we have the
strict inequality $\ov{g} < \ov{g}_{\rm max}$ somewhere on $P^{\circ}$,
it is not too difficult to show that in this case the eigenvalue 
inequalities are strict.

Finally, these results allow us to prove Corollary~\ref{cor1} for metrics
with non-negative Gauss curvature. For such metrics, the corresponding
function $\ov{g}$ satisfies (see \S 2.5)
$$\ov{g}''(x)\leq 0,\ \forall x\in P,$$
and, as an easy consequence, conditions~(\ref{embcond}). Hence, their
invariant eigenvalues are also less than the corresponding ones for
$g_{\rm max}$,
and this is again an optimal bound since the family of ellipsoids
considered above (degenerating to the union of two flat discs) has
positive curvature.

The minimizing property of the diameter of the double flat disc among 
surfaces of revolution in $\R^3$ follows easily 
from~(\ref{eq:diameter}) and~(\ref{embcond}):
$$ D(g) = \int_{-1}^{1} \sqrt{g} dx = \int_{-1}^{1} \frac{\disps 1}
{\disps \sqrt{\ov{g}}} dx > \int_{-1}^{1} \frac{\disps 1}
{\disps \sqrt{\ov{g}_{\rm max}}} = D(g_{\rm max}) = 2 \sqrt{2}\ .$$

\section{Independence of $\la_{1}$ and $D$}

In the case of the metric given in \S~\ref{larginv}, which was used to
obtain a large first invariant eigenvalue, it is an exercise in calculus
to show that
\[
D(g_{\mu}) = {\disps \int_{-1}^{1}} \sqrt{g_{\mu}(x)}dx =
{\disps \int_{-1}^{1}} \frac{\disps dx}{\disps 
\sqrt{(1-x^{2})(1+(1-x^{2})\mu)}}  \to 0
\]
as $\mu$ goes to infinity. Thus, in this case we have that the diameter 
goes to zero while $\lambda_{1}$ goes to infinity. On the other hand, 
it is straightforward to check that in the case where $\lambda_{1}$ 
was made to be arbitrarily small, the corresponding diameter was
unbounded.

As pointed out in the Introduction, these examples together with the fact that the
invariant spectrum is given by a one--dimensional eigenvalue problem, suggest a
natural relation between the diameter of the surface and the value of the first
invariant eigenvalue, in the sense that making the diameter large would give 
small values of $\lambda_{1}$ and vice--versa. However, this is not 
necessarily the case, as we shall now see. In order to do this, we 
need a
result from~\cite{opku} already mentioned in \S~\ref{larginv}. A 
version adequate for our purposes is the following.
\begin{thm} Let $g$ be a continuous even function on $(-1,1)$. The inequality
\[
{\disps \int_{-1}^{1} \frac{\disps 1}{\disps g(x)}
\left[f'(x)\right]^{2}dx} \geq C{\disps \int_{-1}^{1} f^{2}(x)dx}
\]
holds for every odd function $f\in X_{1}$ if and only if
\[
A = \sup_{x\in(0,1)} \left[ (1-x){\disps\int_{0}^{x}}g(t)dt\right]
\]
is finite. Furthermore, the optimal constant $C$ satisfies $A\leq 
1/C\leq 2A$.
\end{thm}
For a more general version of this result, as well as for a proof, 
see~\cite{opku}.

In the case where $g$ is even and corresponds to a metric on the 
sphere, we know that $A$ is finite, and the optimal constant 
corresponds to the first nontrivial eigenvalue, since 
for even $g$'s the corresponding first invariant eigenfunction is odd. 

By appropriate changes of the family of metrics used in 
\S~\ref{larginv}, 
we will now give examples showing that it is possible to have both
the diameter and $\lambda_{1}$ either very small or very large.

\begin{examp}{\rm [$D$ and $\lambda_{1}$ both small] Consider the function
\[
g_{\mu}(x) = \frac{\disps 1}{\disps  (1-x^{2})[1+(1-x^{2})\mu]} +
\frac{\disps \mu^{1-\alpha}}{\disps (1+\mu x^{2})^{2}},
\]
where $\alpha$ is a constant in $(0,1/2)$.
In this case we have
\[
\begin{array}{lcl}
D(g_{\mu}) & = & 2{\disps \int_{0}^{1}} \left[ \frac{\disps 1}{\disps 
(1-x^{2})(1+(1-x^{2})\mu)} + \frac{\disps \mu^{1-\alpha}}{\disps (1+\mu 
x^{2})^{2}} \right]^{1/2} dx\eqskip\\
 & \leq & 2 {\disps \int_{0}^{1}} \frac{\disps 1}{\disps 
\sqrt{(1-x^{2})(1+(1-x^{2})\mu)}} + \frac{\disps 
\mu^{\frac{1-\alpha}{2}}}{\disps 1+\mu 
x^{2}} dx \ \rightarrow 0
\end{array}
\]
as $\mu\rightarrow\infty$ (using the fact that $\alpha > 0$).
Thus we have that the diameter 
converges to zero.

On the other hand, the value of $A$ is now given by 
\[
\begin{array}{lcl}
A & = & {\disps \sup_{x\in(0,1)}} \left\{
(1-x){\disps \int_{0}^{x}} \left[ 
\frac{\disps 1}{\disps  (1-t^{2})(1+(1-t^{2})\mu)}
+ \frac{\disps \mu^{1-\alpha}}{\disps (1+\mu t^{2})^{2}}
\right] \, dt \right\} \eqskip\\
 & \geq &
{\disps \sup_{x\in(0,1)}} \left\{
\mu^{1-\alpha}(1-x){\disps \int_{0}^{x}}\frac{\disps 1}
{\disps (1+\mu t^{2})^{2}}\, dt  \right\} \eqskip\\
 & = &
{\disps \sup_{x\in(0,1)}} \left\{ 
\frac{\disps 1-x}{\disps 2}\left[
\frac{\disps \mu^{1-\alpha} x}{\disps 1+\mu x^{2}} + \mu^{\frac{1}{2}-\alpha}
\arctan\left(\sqrt{\mu}x\right)\right] \right\} \eqskip\\
& \geq &
\frac{\disps 1}{\disps 4}\mu^{\frac{1}{2} - \alpha}
\arctan\left(\sqrt{\mu}/2 \right)
\end{array}
\]
which goes to infinity as $\mu$ goes to infinity, provided $\alpha<1/2$. 
Since $\lambda_{1}\leq 1/A$, we have that
$\lambda_{1}$ goes to zero.
}
\end{examp}

\begin{examp}{\rm [$D$ and $\lambda_{1}$ both large] Consider the function
\[
g_{\mu}(x) = \frac{\disps 1}{\disps  (1-x^{2})(1+(1-x^{2})\mu)} +
\frac{\disps 1}{\disps \log{\mu}}\left( \frac{\disps 1}{\disps 
\left( 1+1/\mu \right)^{2}-x^{2}}\right)^{2}\ .
\]
In this case,
\[
\begin{array}{lcl}
D(g_{\mu}) & = & 2 {\disps \int_{0}^{1}} \left[ \frac{\disps 1}{\disps 
(1-x^{2})[1+(1-x^{2})\mu]} + \frac{\disps 1}{\disps \log{\mu}}\left(
\frac{\disps 1}{\disps
\left( 1+1/\mu \right)^{2}-x^{2}}\right)^{2}
\right]^{1/2}dx\eqskip\\
 & \geq & \frac{\disps 2}{\disps \log^{1/2}{\mu}}
 {\disps \int_{0}^{1}}\frac{\disps 1}{\disps 
 \left[(1+1/\mu)^{2}-x^{2}\right]}dx\eqskip\\
 & = & \frac{\disps \mu}{\disps 1+\mu}\left[\frac{\disps \log(2+1/\mu)}
 {\disps \log^{1/2}\mu} + 
 \log^{1/2}\mu\right]
\end{array}
\]
which goes to infinity as $\mu$ goes to infinity.

For $A$ we now have
\[
\begin{array}{lcl}
A & = & {\disps \sup_{x\in(0,1)}} \left\{(1-x){\disps \int_{0}^{x}}
\left[
\frac{\disps 1}
{\disps  (1-t^{2})(1+(1-t^{2})\mu)}
+ \frac{\disps 1}{\disps \log{\mu}}\left(\frac{\disps 1}{\disps
\left( 1+1/\mu \right)^{2}-t^{2}}\right)^{2} \right]\,dt\right\}\eqskip\\
 & \leq & {\disps \sup_{x\in(0,1)}} \left\{(1-x){\disps \int_{0}^{x}}
 \frac{\disps 1}
{\disps  (1-t^{2})(1+(1-t^{2})\mu)}dt\right\} +\\
& & \hspace*{2cm} + {\disps \sup_{x\in(0,1)}}
\left\{ \frac{\disps (1-x)}{\disps \log\mu}{\disps \int_{0}^{x}}
\frac{\disps 1}{\disps \left(
\left( 1+1/\mu \right) - t\right)^{2} }\,dt\right\}\ \rightarrow 0
\end{array}
\]
as $\mu\rightarrow\infty$, by a simple calculation (note that the first
term corresponds to the family used in \S~\ref{larginv}).
Since $\la_{1}\geq 1/2A$ we conclude that $\la_{1}$ still goes to 
infinity for this family of metrics.
}
\end{examp}

\appendix

\section{Solutions of the Euler-Lagrange equations } 

In the case where $p$ is equal to $1$ (in Lemma~\ref{lemineq}), the second order
differential equation associated with the minimization problem
under study is given by
\begin{equation}
\label{eq2}
-\frac{\disps d}{\disps dx}\left[ (1-x^2)^2
\frac{\disps df}{\disps dx}\right] = \lambda f.
\end{equation}
It turns out that this equation can actually be solved explicitly
for all real values of the parameter $\lambda$. To see this, we
begin by writing it as
\[
(1-x^{2})^{2}\frac{\disps d^{2}f}{\disps dx^{2}}-4x(1-x^{2})
\frac{\disps df}{\disps dx} + \lambda f=0, \ \ x\in(-1,1).
\]
Following~\cite{cole}, we now look for solutions of the form
\[
f(x)=(1+x)^{r}\left[\alpha_{0}+\alpha_{1}(1+x)+\ldots\right]
\]
from which we obtain the indicial equation
\[
r^{2}+r+\lambda/4=0
\]
with two solutions
\[
r_{\pm} = \frac{\disps -1\pm i\omega}{\disps 2}, \ \
\omega=\sqrt{\lambda-1}.
\]
We then see that
\[
f(x) = (x+\omega)(1-x)^{r_{-}}(1+x)^{r_{+}}
\]
satisfies~(\ref{eq2}). Due to the symmetry of the problem, $f(-x)$ is
also a solution. The form of the (real) solutions will now depend on
the sign of $1-\lambda$.

In the case where $\lambda>1$, we have that $r_{\pm} = (-1\pm wi)/2$
and it is possible to obtain from the general form of the solution
that:
\[
f_{1}(x) = \frac{\disps 1}{\disps \sqrt{1-x^{2}} }\left[x\cos\left[
\frac{\disps \omega}{\disps 2}\log\left(\frac{\disps 1-x}{\disps
1+x}\right)
\right] - \omega\sin\left[ \frac{\disps
\omega}{\disps 2}\log\left(\frac{\disps 1-x}{\disps 1+x}\right)
\right]\right]
\]
and
\[
f_{2}(x) = \frac{\disps 1}{\disps \sqrt{1-x^{2}}
}\left[\omega\cos\left[ \frac{\disps \omega}{\disps
2}\log\left(\frac{\disps 1-x}{\disps 1+x}\right)
\right] + x\sin\left[ \frac{\disps
\omega}{\disps 2}\log\left(\frac{\disps 1-x}{\disps 1+x}\right)
\right]\right].
\]
are two linearly independent real solutions (odd and even,
respectively).

When $\lambda=1$ we get that
\[
\begin{array}{lcl}
f_{1} = \frac{\disps x}{\disps \sqrt{1-x^{2}} } & \mbox{and} & f_{2} =
\frac{\disps 1}{\disps \sqrt{1-x^{2}} }\left[ \frac{\disps x}{\disps
2}\log\left(\frac{\disps 1-x}{\disps 1+x}\right)+1\right]
\end{array}
\]
are odd and even solutions, respectively.

When $\lambda<1$, and letting now
$\gamma=-i\omega=\sqrt{1-\lambda}$, we obtain the two solutions
\[
f_{1}(x) = (\disps x+\gamma)(1+x)^{r_{-}} (1-x)^{r_{+}}+
(x-\gamma)(1+x)^{r_{+}} (1-x)^{r_{-}}
\]
and
\[
f_{2}(x) = (\disps x+\gamma)(1+x)^{r_{-}} (1-x)^{r_{+}}-
(x-\gamma)(1+x)^{r_{+}} (1-x)^{r_{-}}.
\]

>From this we see that there are no nontrivial solutions of the
Euler--Lagrange equations lying in the space $X$.





\begin{thebibliography}{9999}

\bibitem[A]{abr} M. Abreu, `K\"ahler geometry of toric manifolds in 
symplectic coordinates', preprint (2000), math.DG/0004122.

\bibitem[B]{berg} M. Berger, `Encounter with a geometer: Eugenio 
Calabi', in {\em Manifolds and Geometry}, 20--60 (eds. P. de Bartolomeis, F. 
Tricerri and E. Visentini), Symposia Mathematica (Vol.XXXVI), 
Cambridge University Press 1996.

\bibitem[BLY]{bly} J.-P. Bourguignon, P. Li and S.T. Yau, `Upper 
bound for the first eigenvalue of algebraic submanifolds', Comment. 
Math. Helvetici {\bf 69} (1994), 199--207.

\bibitem[CL]{cole} E. Coddington and N. Levinson, {\em Theory of Ordinary
Differential Equations}, McGraw--Hill 1955.

\bibitem[E1]{eng1} M. Engman, `Trace formulae for $S^1$ invariant 
Green's operators on $S^2$', Manuscripta Math. {\bf 93} (1997), 357--368.

\bibitem[E2]{eng2} M. Engman, `The spectrum and isometric embeddings of 
surfaces of revolution', preprint (1999), math.DG/9910038.

\bibitem[G]{guil} V. Guillemin, `K\"{a}lher structures on toric
varieties', J. Differential Geometry {\bf 40} (1994), 285--309.

\bibitem[H]{hers} J. Hersch, `Quatre propri\'{e}t\'{e}s
isop\'{e}rim\'{e}triques de membranes sph\'{e}riques homog\`{e}nes', C.R.
Acad. Sci.Paris S\'{e}r. A-B {\bf 270} (1970), A1645--A1648.

\bibitem[KW]{kawa} J. Kazdan and F. Warner, `Surfaces of revolution
with monotonic increasing curvature and an application to the equation
$\De u = 1-Ke^{2u}$ on $S^2$', Proc. Amer. Math. Soc. {\bf 32} (1972), 139--141.

\bibitem[N]{nire} L. Nirenberg, `The Weyl and Minkowski problems in 
differential geometry in the large', Comm. Pure Appl. Math. {\bf 6} 
(1953), 337--394.

\bibitem[O]{onei} B. O'Neill, {\em Elementary differential geometry},
Academic Press 1966.

\bibitem[OK]{opku} B. Opic and A. Kufner, {\em Hardy--type 
inequalities},
Pitman Research Notes in Mathematic Series {219}, Longman Scientific
and Technical, Harlow 1990.

\bibitem[YY]{yaya} P. Yang and S.T. Yau, `Eigenvalues of the Laplacian 
of Compact Riemann Surfaces and Minimal Submanifolds', Ann. Scuola Norm. 
Sup. Pisa Cl. Sci. {\bf 7} (1980), 55--63.

\bibitem[Y]{yau} S.T. Yau, {\em Seminar on differential geometry}, Ann. Math. 
Stud. {\bf 102}, Princeton University Press, Princeton, N.J., 1982.

\end{thebibliography}
\end{document}